\documentstyle[12pt,epsf]{article}

\font\twelve=cmbx10 at 15pt
\font\ten=cmbx10 at 12pt

\renewcommand{\thefootnote}{\fnsymbol{footnote}}
\newcommand\be{\begin{equation}}
\newcommand\ee{\end{equation}}
\newcommand{\ba}{\begin{array}}
\newcommand{\ea}{\end{array}}
\newtheorem{ttheorem}{{\bf Theorem}}[section]
\newtheorem{llemma}[ttheorem]{{\bf Lemma}}
\newtheorem{prop}[ttheorem]{{\bf Proposition}}
\newtheorem{cor}[ttheorem]{{\bf Corollary}}
\newenvironment{corollary}{\begin{cor}\begin{rm}}{\end{rm}\end{cor}}
\newenvironment{theorem}{\begin{ttheorem}\begin{rm}}{\end{rm}%
\end{ttheorem}}
\newenvironment{lemma}{\begin{llemma}\begin{rm}}{\end{rm}\end{llemma}}
\newenvironment{proposition}{\begin{prop}\begin{rm}}{\end{rm}\end{prop}}
\def\wt{\widetilde}
\def\wh{\widehat}

\newcommand{\rf}[1]{\ref{#1}]}

\begin{document}

\begin{titlepage}

\begin{center}

{\ten Centre de Physique Th\'eorique\footnote{Unit\'e Propre de
Recherche 7061} - CNRS - Luminy, Case 907}

{\ten F-13288 Marseille Cedex 9 - France }

\vspace{1,5 cm}

{\twelve A NEW APPROACH TO CHARACTERISE \\
ALL THE TRANSITIVE ORIENTATIONS \\
FOR AN UNDIRECTED GRAPH}

\vspace{0.3 cm}

\setcounter{footnote}{0}

\renewcommand{\thefootnote}{\arabic{footnote}}

{\bf Melhem HAMAD\'E\footnote{Universit\'e de Provence} and
Abdelkhalek BELKASRI\footnote{and Universit\'e d'Aix-Marseille II} }

\vspace{2 cm}

{\bf Abstract}

\end{center}

A new approach to find all the transitive orientations for a
comparability graph (finite or infinite) is presented.
This approach is based on the link between the notion
of ``strong'' partitive set and the forcing theory
(notions of simplices and multiplices). A mathematical algorithm is
given  for the case of a comparability graph which has only non limit
sub-graphs.

\vspace{2 cm}


\bigskip

\noindent Number of figures : 7

\bigskip

\noindent November 1994

\noindent CPT-94/P.3069

\bigskip

\noindent anonymous ftp or gopher: cpt.univ-mrs.fr

\end{titlepage}

\section{Introduction}
 The problem of the transitive orientation for a comparability graph
 was studied by Golumbic
using the forcing theory [1]. The problem was solved for {\it finite}
comparability graphs and an algorithm was given which gives one
transitive orientation for a {\it finite} comparability graph.

The purpose of this paper is to study the transitive orientation for
the case of {\it infinite} comparability graphs. The results for the
finite case could not be extended to the infinite case because of the
finite type like of the approach used in Ref.[1]. We were then
obliged to consider a new approach, but remaining within the forcing
theory. This was possible by introducing the notion of a {\it
``strong'' partitive set}. It happens that this idea permits to solve
the problem of the transitive orientations by inducing a lot of
characteristics of the forcing theory, from one hand and the
undirected graphs, in general, from the other hand.

The paper is organized as follows: the section 2 is devoted to the
definitions and notations used throughout the paper and also to recall
the main results of the forcing theory which still valid for infinite
graphs. In section 3, we establish narrow links between the notion of
a {\it ``strong'' partitive set} and {\it simplices} (and {\it
multiplices}) which are the principal touls of the forcing theory.
In section 4, we used the results of section 3 to face the problem of
transitive orientation of comparability graphs. We proved a theorem
which is in fact  a mathematical algorithm which gives all the
transitive orientations for a comparability graph
which has all its sub-graphs non limit.

\section{Preliminairies}
This section is devoted to the definitions and notations which will be
used throughout this article. We also recall some results about the
partitive sets and the implication classes. These results can be found
with more details in Ref.[1].

We consider here any kind of graphs; finite or infinite. In what
follows we denote by $G=(V,E)$ any graph, where $V$ is the set of
vertices, and by $E(\subseteq V^2)$ the set of edges. Directed edge
will be denoted by $(a,b)$ \newline
(for $a,b\in V$) and an undirected one is denoted
by $ab, ab=\left\{(a,b),(b,a)\right\}$. We say that a graph $G=(V,E)$
is {\it empty} if $E=\phi$.

For any $X\subseteq V$, $G(X)=(X,E(X))$ will denote the sub-graph of
$G$ induced on $X$; where
$$
E(X)=\left\{(a,b)\in E\quad ;\quad \{a,b\}\subseteq X\right\}
$$
we define also the set of vertices $\widetilde{A}$ spanned by a set of
edges $A$ as
$$
\widetilde{A} =\left\{a\in V;\ \hbox{there exists}\ b\in V
\hbox{so that}\ (a,b)\in A
\
\hbox{or}\ (b,a)\in A\right\}
$$

\subsection{Partitive set and ``strong'' partitive set}
Let $\cong$ be a binary relation acting on $V^2$ defined by
$$(a,b)\cong (c,d)\Leftrightarrow\left(\{(a,b),(c,d)\}\cap E=\phi\
\hbox{or}\ \{(a,b),(c,d)\}\subseteq E\right)
$$
This means that the edges $(a,b)$ and $(c,d)$ are both belonging to
$E$ or are both out of $E$.

Let $X\subseteq V$ be a sub-set of $V$. $X$ is a partitive set of $G$
(or $V$) if:
$$
\hbox{for every}\ \{a,b\}\subseteq X\ \hbox{and every}\ c\in V-X
$$
$$
\hbox{we have}\ (a,c)\cong (b,c)\ \hbox{and}\ (c,a)\cong (c,b)
$$
It means that the elements of $X$ are related to any external element
in the same manear. The notion of the partitive set is the analogue of
the notion of interval in an ordered set. We will denote by $I(G)$ the
class of partitive sets of the graph $G$. A partitive set is trivial
if it is a singleton or equal to $V$. By $I^{\star}(G)$ we will denote
the class of  non-trivial partitive sets.

\noindent
A graph is indecomposable if all its  partitive sets are
trivial; otherwise, it is decomposable.

\noindent
A partitive set $X\in I(G)$ is called a ``strong'' partitive set of
$G$(or $V$) if
for every partitive set $Y\in I(G)$ so that $Y\cap X\not =\phi$,
we have either $X\subseteq Y$ or $Y\subseteq X$.
It means that the eventuality that $X$ and $Y$ shear only a common
part is excluded.

\noindent
We will denote by $I_F(G)$ the class of ``strong'' partitive sets of
$G$ and by $I^{\star}_F(G)$ the class of the non-trivial ``strong''
partitive sets. We say that $G$ is  {\it limit} if $I_F(G)$ do not
contain any element different from $V$
 which is maximal for the inclusion$(\subseteq)$; otherwise it is
{\it non limit}.

\newpage
\noindent {\it Isomorphism}: Two graphs $G=(V,E)$ and $G'=(V',E')$ are
said to be isomorphic if:

\begin{enumerate}

\item[(i)] {there is a bijection $f:V\to V'$;}
\item[(ii)] {$f$ preserves the edges, i.e., for every
$\{a,b\}\subseteq V:$}
\end{enumerate}
$$
(a,b)\in E\Leftrightarrow (f(a),f(b))\in E'\ .
$$

\noindent {\it Quotient graph}: Let $G=(V,E)$ be a graph and $P$ a
partition of $G$ made of partitive sets $(P\subseteq I(G))$. We define
the {\it quotient graph} of $G$ by $P$, denoted by $G/P$, as the
isomorphic graph to $G(f(P))$, where $f$ is a choice function from $P$
to $V$, i.e., $X\in P\Rightarrow f(X)\in X$.

\begin{proposition} \label{prop2.1}
 Let $P$ be a partition of partitive sets (respectively of ``strong''
partitive sets) of $G=(V,E)$ and $X$ a sub-set of $P$. We have then\\
$X\in I(G/P)$ (respectively $X\in I_F(G/P)$) if and only if:
$\cup X\in I(G)$

\noindent $(\hbox{respectively}\ \cup X\in I_F(G))\ ,$

\noindent
where $\cup X$ means the union of the vertices constituting the
partitive sets (respectively the ``strong'' partitive sets) of $X$.
\end{proposition}

\subsection{Implication classes and simplices}
{\it Comparability graph} {\it (or transitively orientable graph)}:
Let
$G=(V,E)$ be an undirected graph. $G$ is a comparability graph if
there exists an orientation of the edges of $G$ so it constitutes a
partial order on $V$.

Comparability graphs are also known as transitively orientable graphs
or partially orderable graphs.

\noindent{\it Implication classes}: Let us define the binary relation
$\Gamma$ on the edges of an undirected graph $G=(V,E)$ as follows:
$$
(a,b)\Gamma(a',b')\Leftrightarrow\left\{\begin{array}{llll}
\hbox{either}& a=a'&\hbox{and}&bb'\not\in E\\
\hbox{or}&b=b'&\hbox{and}&aa'\not\in E.\\
\end{array}\right.
$$
We say that $(a,b)$ directly forces $(a',b')$ whenever $(a,b)\Gamma
(a',b')$. In graphical representation $ab$ and $a'b'$ will have only
one common vertex and the orientation is so that arrows are both
pointed to the extremities or are both pointed to the common vertex
(see Figure 2.1).

Notice that $\Gamma$ is not transitive. In Figure~2.1 $(a,d)$ forces
directly $(a,b)$ and $(a,c)$. The edges $(b,c)$  and $(c,b)$ are
 not forced by any
other edge.

\begin{center}

\vglue 0,2cm


\vglue 0,4cm

{\bf Figure 2.1}

\end{center}

Let $\Gamma^{\star}$ be the transitive closure of $\Gamma$. It is easy
to verify that $\Gamma^{\star}$ is an equivalence relation on $E$ and
the equivalence classes related to $\Gamma^{\star}$ are what one
call the {\it implication classes} of $G$.

In what follows we will see that it is useful to define what one call
color classes of $G$ (or shortly colors of $G$). If $A$ is an
implication class, the color class associated to $A$ and denoted by
$\widehat A$, is the union of $A$ and $A^{-1}$ \\
$(\widehat A=A\cup
A^{-1}\subseteq E)$; where $A^{-1}=\{(a,b)\in E$ with $(b,a) \in A\}$.

\begin{theorem}\label{th2.2}
({\it Golumbic Ref.[1]})  Let $G=(V,E)$ be a comparability graph and
$A$ an implication class of $G$. If $O=(V,E')$ is a partial order
associated to $G$, we have necessarily either $E'\cap\widehat A=A$ or
$E'\cap
\widehat A=A^{-1}$ and, in either case, $A\cap A^{-1}=\phi$.
\end{theorem}

\begin{lemma}\label{lem2.3} {\it (The Triangle Lemma).} Let
$A,B$ and
$C$ be implication classes of an undirected graph $G=(V,E)$ with $A
\not = B$ and $A\not = C^{-1}$ and having edges $(a,b)\in C$,
$(a,c)\in B$ and
$(b,c)\in A$, we have then

\begin{enumerate}

\item[(i)] $(b',c')\in A\Rightarrow \left((a,b')\in C\right.$ and
$\left.(a,c')\in B\right)$ ;

\item[(ii)] $\left((b',c')\in A\right.$ and $\left.(a',b')\in
C\right)\Rightarrow (a',c')\in B$ ;

\item[(iii)] $a\not\in\widetilde A$.

\end{enumerate}

\end{lemma}

The following results are consequences of the triangle lemma.

\begin{center}

\vglue 0,2cm


\vglue 0,4cm

{\bf Figure 2.2 } {\it Illustration of the triangle lemma.}
\end{center}

\begin{theorem}\label{th2.4}
({\it Golumbic Ref.[1]})  Let $A$ be an implication class of an
undirected graph $G=(V,E)$. Exactly one of the following alternatives
holds:
\begin{enumerate}
\item[(i)] $A=A^{-1}=\widehat A$ and $\widehat A$ is not transitively
orientable ;
\item[(ii)] $A\cap A^{-1}=\phi$ and then $A$ and $A^{-1}$ are two
transitive orientations of $\widehat A$.
\end{enumerate}

\end{theorem}

\begin{proposition}\label{prop2.5}
Let $X$ be a partitive set $(X\in I(G))$ and $\widehat A$ a color
class of an undirected graph $G=(V,E)$ so that $E(X)\cap\widehat
A\not=
\phi$ we have then
$\widehat A\subseteq E(X)$.
\end{proposition}

\begin{proposition}\label{prop2.6}
If $\widehat A$ is a color class of $G=(V,E)$ then $\widetilde A$ is a
partitive set of $G$ $(\widetilde A\in I(G))$.
\end{proposition}

\noindent {\it Simplex} : Let $G=(V,E)$ be an undirected graph. A
$K_{r+1}$ complete sub-graph of $G$, $S=(V_S,E_S)$ on $r+1$ vertices
is called a {\it simplex of rank} $r$ if each undirected edge $ab$ of
$E_S$ is contained in a different color class of $G$. A simplex is
{\it maximal} if it is not properly contained in any larger simplex.

The {\it multiplex} $M$ generated by a simplex $S$ of rank $r$ is
defined to be the part of $E$ constituted of all edges which their
color classes are present in the simplex $S$
($M(S)={{\cup}_{\widehat A\cap E_S \neq\phi}}\widehat A$).
$M$ is said also a
multiplex of rank $r$. The multiplex $M$ is said to be maximal if $S$
is maximal. We will denote by $\widehat M$ the collection of color
classes present in the multiplex $M$.

\section{Connection between multiplices and ``strong'' partitive sets}
In this section we make connection between the notion of ``strong''
partitive set and multiplices. In our knowledge this was never made
before. It happens that this connection allows us to recover results
of Golumbic [1] for finite graphs and generalize them to the infinite
graphs. The results presented in this section will be used in the
following section to state a theorem on the decomposability of
undirected graphs. In what follows $G$ will denote an undirected graph
unless other mention is pointed out.

\begin{proposition}\label{prop3.1}
If $\widehat A$ and $\widehat B$ are two color classes of $G=(V,E)$ we
have then $(\widehat A=\widehat B)\Leftrightarrow (\widetilde
A=\widetilde B)$.

\end{proposition}

{\it Proof.} It is obvious that $\widehat A=\widehat B \Rightarrow
\widetilde A=\widetilde B$. Let us suppose that $\widetilde A=
\widetilde B$ and $\widehat A\not=\widehat B$, we have then for every
$x\in\widetilde A$ there exists $a\in\widetilde A$ and $b\in\widetilde
B$ so that $ax\in \widehat A$ and $bx\in\widehat B$. Since $\widehat
A\not=\widehat B$ we have necessarily $ab\in E$. Let $C$ be the color
class which contains $ab$. We have then two alternatives:

\begin{enumerate}

\item[(i)] $\wh C\not=\wh A\not= \wh B\Rightarrow b\not\in\wt A$
(Triangle Lemma) which is absurd since\\
 $b\in \wt B=\wt A$;

\item[(ii)] $\wh C=\wh A\not= \wh B\Rightarrow a\not\in\wt B$
(Triangle Lemma) which is also absurd, since $a\in \wt A=\wt B$.

\end{enumerate}

\begin{proposition}\label{prop3.2} Let $\wh A$ be a color class of
$G=(V,E)$ and $X$ a partitive set of $G$ so that $X\subset \wt A$,
there exists then $a\in\wt A-X$ so that for every $x\in X$, $ax\in\wh
A$.

\end{proposition}

{\it Proof.} Let $x\in X\subset\wt A$, there exists $a\in\wt A$ so
that $ax \in \wh A$. Necessarily $a\in \wt A-X$ (otherwise $(ax\in\wh
A\cap E(X))\Rightarrow\wh A\subseteq E(X)$~[P.\rf{prop2.5}) then for
every
$x\in X$ we have $xa\in E$, since $X\in I(G)$. Let us assume that
there exists $y\in X$ and a color $ \wh B\not=\wh A$ so that $ay\in
\wh B$, then $xy\in E$. Let $\wh C$ be the color which contains $xy$.
We have
$\wh C\not=\wh A$ and $\wh C\not=\wh B$ (otherwise: $(\wh C=\wh
A\Rightarrow \wt A\subseteq X)$ and $(\wh C=\wh
B\Rightarrow \wt B\subseteq X)$ which contradicts the fact that
$a\in(\wt A\cap\wt B)-X$).
Then using the triangle lemma we will have $y\not\in\wt A$ which is
absurd. Finally we have for every $x\in X,ax\in \wh A$.

\begin{theorem}\label{th3.3}
Let $\wh A$ and $\wh B$ be two color classes of $G=(V,E)$ so that $\wt
A-\wt B\not=\phi$ and $\wt B-\wt A\not=\phi$. We have then that $X=\wt
A\cap \wt B$ is a ``strong'' partitive set of $G$.

\end{theorem}

{\it Proof.} Since $\wt A$ and $\wt B$ are partitive sets[P.2.6] we
have that
$X=\wt A\cap \wt B$ is a partitive set. If $X=\phi$ or $X$ is a
singleton $(|X|=1)$ then $X$ is a ``strong'' partitive set.

Let us suppose now that $|X|>1$. Let $Y$ be a partitive set $(Y\in
I(G))$ so that $X\cap Y\not=\phi$ and $Y-X\not=\phi$ and let
$z\in X\cap Y$ and $y\in Y-X$. We have to show
that $X\subset Y$.
Applying~[P.\rf{prop3.2} we have:

$X\subset\wt A\Rightarrow$ there exists $a\in\wt A-X$ so that for
every
$x\in X,\ ax\in\wh A$;

$X\subset\wt B\Rightarrow$ there exists $b\in\wt B-X$ so that for
every
$x\in X,\ bx\in\wh B$.

Thus $za\in\wh A$ and $zb\in \wh B$.

If $a\in Y$ then $az\in E(Y)\cap\wh A\not=\phi$, thus $X\subseteq\wt
A\subseteq Y$~[P.\rf{prop2.5}.

If $b\in Y$ then $bz\in E(Y)\cap\wh B\not=\phi$, thus $X\subseteq\wt
B\subseteq Y$~[P.\rf{prop2.5}.

Let us now suppose that $\{a,b\}\cap Y=\phi$, we have then:

\begin{enumerate}

\item[(i)] if $y\not\in\wt A$ then: $(\{y,z\}\subseteq Y\in I(G)$ and
$za\in E)\Rightarrow ya\in E$; \\
$(\wt A\in I(G)$ and $ay\in E)\Rightarrow$
for every $x\in X\subset \wt A, yx\in E$. \\
Let $\wh K$ be the color class of $ya$ and $\wh C$ that one of
$yz$. We have then: $(y\in(\wt K\cap\wt C)-\wt A$, $a\in\wt K-Y$ and
$yz\in\wh C\cap E(Y)\not=\phi)$ implies that  $(\wt K\not=\wt A\not=
\wt C)\Leftrightarrow (\wh K\not=\wh A\not=\wh C)$, where at the last
step we have applied~[P.\rf{prop3.1}. \\
Let $v\in X$ (with $v\not= z$) and $\wh D$ be the color of $yv$. We
have then: \\
$\wh K\not=\wh A\not=\wh D$ (since $\wh K\not=\wh A,yv\in\wh D\cap
E(Y\cup X)\not=\phi$ and\\
 $a\in (\wt A\cap\wt K)-(Y\cup X)$).
The two tricolor triangles $(a,y,z;\wh K,\wh A,\wh C)$ and $(a,y,v;\wh
K,\wh A,\wh D)$ have two common colors, thus $\wh D=\wh C$, thus for
every $x\in X, yx\in\wh C$ which implies that $X\subset\wt C\subseteq
Y$ ($\wt C\subseteq Y$ because $yz\in \wh C\cap E(Y)\neq\phi$[P.2.5])
 (see Figure~3.1).

\item[(ii)] If $y\in \wt A$ then $(y\not\in X=\wt A\cap\wt B$ and
$y\in\wt A)\Rightarrow y\not\in\wt B$. By replacing $a$ by $b$ and
$\wh A$ by $\wh B$ up here in (i), we get $X\subset \wt C \subseteq Y$.

\end{enumerate}

\begin{center}

\vglue 0,2cm


\vglue 0,4cm

{\bf Figure 3.1} {\em We have : $\wh D=\wh C$[Triangle Lemma].}

\end{center}

\smallskip

\begin{theorem}\label{th3.4}
Let $X$ and $Y$ be two ``strong'' partitive sets of $G=(V,E)$ so that
$X\cap Y=\phi$. $X$ and $Y$ can be related to each other by only one
color at most.

\end{theorem}

\begin{center}

\vglue 0,2cm


\vglue 0,4cm

{\bf Figure 3.2}
{\em If $X$ and $Y$ are ``strong'' partitive sets we have necessarily
$\wh A = \wh B$}.

\end{center}

\medskip

{\it Proof.} Let $\wh A$ and $\wh B$ be two colors connecting $X$ and
$Y$. We have then $X\subset \wt A\cap\wt B$ and $Y\subset \wt A\cap
\wt B$ (applying [P.2.6] and  the definition of the ``strong''
partitive set). Let us suppose that $\wh A\not=\wh B$ and define the
following sets
$$I_A=\left\{Z\in I_F(G)\ \hbox{with}\ Z\subset\wt A\right\}\
\hbox{and}\ I_B=\left\{Z\in I_F(G)\ \hbox{with}\ Z\subset\wt
B\right\}$$ we have $\cup I_A=\wt A$ and $\cup I_B=\wt B$ since:
$\cup I_A\subseteq \wt A$ and
if $x\in \wt A$, then $\{x\}\in I_F(G)\cap I_A$, thus $x\in \cup I_A$
and then $\wt A\subseteq \cup I_A$. The same situation holds for $I_B$.

Since $\{X,Y\}\subseteq I_A\cap I_B$ we have $\wt A\cap \wt
B\not=\phi$, $\wh A\cap E(\wt A\cap \wt B)\not=\phi$ and $\wh B\cap
E(\wt A\cap \wt B)\not= \phi$ we have then only two alternatives:

\begin{enumerate}
\item[(i)] either $\wt A=\wt B\Rightarrow\wh A=\wh B$~[P.\rf{prop3.1}
\item[(ii)] or $\wt A\not=\wt B\Rightarrow\left(\wt A-(\wt A\cap \wt
B)\neq\phi
\right.$ or $\left.\wt B -(\wt A\cap \wt B)\neq\phi\right)$ which is
absurd since it contradicts with $\wh A\cap E(\wt A\cap\wt B)\not=
\phi$ and $\wh B\cap E(\wt A\cap\wt B)\not=\phi$~[P.\rf{prop2.5}.

\end{enumerate}

\begin{corollary}\label{cor3.5}
Let $P$ be a partition of ``strong'' partitive sets of $G=(V,E)$ and
$f$ a choice mapping from $P$ to $V$, i.e., $f:X\in P\to f(X)\in
X\subseteq V$. We have then that the isomorphism from $G/P$ to
$G(f(P))$ conserves the color classes.
\end{corollary}

\begin{corollary}\label{cor3.6}
Let $X$ be a ``strong'' partitive set of $G=(V,E)$ with $X\not= V$ and
let $u\in V-X$. We have then:
\begin{enumerate}
\item[(i)] either for every $x\in X, ux\not\in E$;
\item[(ii)] or there exists a color $\wh A$ of $G$ so that for every
$x\in X, ux\in \wh A$.
\end{enumerate}
\end{corollary}

{\it Proof.} This is true because for every $u\in V$ the singleton
$\{u\}$ is a ``strong'' partitive set of $G$.

\begin{theorem}\label{th3.7}
Let $X$ be a ``strong'' partitive set of $G=(V,E)$ and $M(S)$ a
multiplex of $G$ generated by $S=(V_S,E_S)$. We have then the
following implication:
$M\cap E(X)\not=\phi\Rightarrow M\subseteq E(X)$.

\end{theorem}

{\it Proof.} If $E_S\subseteq E(X)$ then for every color $\wh
A\subseteq M$ we have $\wh A\cap E(X)\not =\phi$ and thus for every
$\wh A\subseteq M$, $\wh A\subseteq E(X)$ which implies that
$M\subseteq E(X)$.

In the other hand if $E_S\cap E(X)=\phi$ then $M\cap E(X)=\phi$. Let
us assume now that $E_S-E(X)\not=\phi$ and $E_S\cap E(X)\not=\phi$,
then $V_S-X\not=\phi$. Let $u\in V_S-X$ and $ab\in E_S\cap E(X)$,
then $\{a,b\}\subseteq X\cap V_S$ and $u$ is related to $a$ and $b$ by
two different colors
(applying the definition of a simplex). This is absurd because it
contradicts with the corollary~[\rf{cor3.6}.

This result is the analogue for multiplices and "strong" partitive
sets  of [P.2.5] which deals with colors and partitive sets.

\begin{lemma}\label{lem3.8}
Let $G=(V,E)$ be an undirected graph with the number of vertices
greater than 2 $(|V|>2)$. If $G$ is decomposable and has a color $\wh
A$ so that $\wt A= V$, then $G$ has a non-trivial maximal ``strong''
partitive set.

\end{lemma}

{\it Proof.} Let $X$ be a non-trivial partitive set of $G$ $(X\in
I^{\star}(G))$, thus $X\not=\wt A$. Let us define on $I(G)$ the
following binary relation $R$  : \newline
$XRY\Leftrightarrow ((X=Y)$ or $(X\cap
Y\not=\phi, X-Y\not=\phi$ and $Y-X\not=\phi))$. And let $R^{\star}$ be
the transitive closure of $R$.

It is easy to verify that $R^{\star}$ is an equivalence relation on
$I(G)$. Let $X^{\star}$ be the equivalence class of $X$ modulo
$R^{\star}$. By the definition itself of $R^{\star}$, $\cup X^{\star}$ is a
``strong'' partitive set of $G$.

In the other hand $X\not=\wt A\Rightarrow$ for every $Y\in X^{\star}$,
$E(Y)\cap\wh A=\phi$ (otherwise $\wt A=Y$ and for every $Z\in
X^{\star}$, $Z\subset Y$ which contradicts with the definition of
$R^{\star}$). We have then $\cup X^{\star}\subset V$. There exists
then
$a\in V-\cup X^{\star}$ so that for every $x\in \cup X^{\star}$,
$ax\in\wh A$[P.3.2]. But $\{ a\}\in I^{\star}_F(G)$ and if $Y\in
I^{\star}_F(G)$ so that  $\{ a\}\subset Y$ then $Y\subseteq V-\cup
X^{\star}$ ( otherwise $\wh A\cap E(Y)\neq\phi\Rightarrow\wt A
=V\subseteq Y $[P.2.5]), consequently $a$ is contained in a maximal
"strong" partitive set different from $V$.  Thus $\cup X^{\star}$ is
contained in a non trivial maximal "strong" partitive set.

\begin{theorem}\label{th3.9}
Let $M (S)$ be a multiplex of $G = (V,E)$ generated by the simplex $S
= (V_S,E_S)$. $M (S)$ is maximal if and only if the set of vertices
$\wt M$ spanned by $M$ is a ``strong'' partitive set.

\end{theorem}

{\it Proof.} Let us suppose that $M (S)$ is maximal, then $S$ is
maximal. $\wt M$ is a partitive set of $G$. Let $Y$ be a partitive set
of $G$ so that $Y \cap \wt M \not= \phi$ and $Y - \wt M \not= \phi$.
Assume that $\wt M - Y \not= \phi$ and let us show that it is absurd.
Let $y \in Y - \wt M$. $G (\wt M)$ is connected and then there exists
$u
\in Y \cap \wt M$ and $v \in \wt M - Y$ so that $u v \in M$. The
following statements hold:
\arraycolsep2.5pt
$$
\ba{ll}
( Y \in I (G),\ u \in Y\ \mbox{and}\ v \in \wt M - Y ) &
\Rightarrow\ \mbox{For any} x\in Y, xv\in E\Rightarrow\ y v \in E\ ;
\\[3mm]
( {\wt M} \in I (G),\ v \in \wt M\ \mbox{and}\ y \in Y - \wt M ) &
\Rightarrow\ \mbox{for any}\ x \in \wt M,\ x y \in E
\ea
$$(see Figure 3.3).

\begin{center}

\vglue 0,2cm


\vglue 0,4cm

{\bf Figure 3.3}
{\em For every $x\in {\wt M},$ $xy\in E$}.

\end{center}

\medskip

The colors connecting $y$ to the summits of $S$ can not be all
different, otherwise, $M$ will not be maximal.
Let us suppose that there exists $\{ a,b \} \in S$ and a color $\wh A$
so that $\{ ya, yb \} \subseteq \wh A$. Since $y \in \wt A - \wt M$
and
$\wt M \in I (G)$, then we have $\wh A \cap M = \phi$[P.2.5]. If
rank$(M)=1$ then there exists a color $\wh K$ so that $\wh K = M$. But
$(ab\in \wh K
\cap E (\wt A) \not= \phi ) \Rightarrow \wt K \subseteq \wt A$ :
$\{ uv,ab\}\subseteq\wh K \Rightarrow \{ yu,yv\}\in\wh A$[Triangle
Lemma], then we have
$(y u \in \wh A \cap E (Y)) \Rightarrow \wt M = \wt K \subseteq \wt
A \subseteq Y$[P.2.5] ( see Figure 3.4) which is
in contradiction with our proposition.

\begin{center}

\vglue 0,2cm


\vglue 0,4cm

{\bf Figure 3.4}
{\em $yu\in\wh A\cap E(Y)\neq\phi$ and $ab\in\wh K\cap E(\wt A)\neq
\phi$}.

\end{center}

\medskip

Let us assume now that rank$(M) \ge 2$. Let $c \in V_S - \{ a,b \},\
\wh B$ the color of $y c,\ \wh C$ the color of $a c$ and $\wh D$ the
color of $b c$. If $\wh A \not= \wh B$ then $S$ contains two tricolor
triangles: $(y,a,c~;\ a y \in \wh A,\ c y \in \wh B,\ a c \in \wh C)$
and $(y,b,c~;\ b y \in \wh A,\ c y \in \wh B,\ b c \in \wh D)$ having
two common colors, then $\wh C = \wh D$~[Triangle lemma] which is
absurd because $S$ is a simplex, then $\wh A = \wh B$ (see Figure
3.5).

\begin{center}

\vglue 0,2cm


\vglue 0,4cm

{\bf Figure 3.5}
{\em $\wh B\neq\wh A\Rightarrow \wh C=\wh D$}.

\end{center}

\medskip


Thus for every $x \in V_S,\ y x \in \wh A$ and by consequence $V_S
\subset \wt A$  which implies that for every color $\wh H \subset M$
we have $ \wh H \cap E (\wt A) \not= \phi$ then $\wt M \subset \wt A$.
But
$(y u \in \wh A \cap E(Y) \not= \phi) \Rightarrow \wt A \subseteq Y$
hence $\wt M \subset \wt A \subseteq Y$, which also contradicts the
first assumption $(\wt M - Y \not= \phi)$.

Finally we have showed that $\wt M  \subset Y$ which means
that $\wt M$ is a ``strong'' partitive set. Let us now prove the
converse. Let us assume that $\wt M$ is a ``strong'' partitive set of
$G$ and let us show that $M$ is a maximal. If the summits of $\wt M$
are related to
$y \in V - \wt M$, then they are related by the same color which
achieves the proof.

\begin{corollary}\label{cor3.10}
Let $M (S)$ be a multiplex of $G = (V,E)$ generated by the simplex
$S=(V_S,E_S)$ and let $a \in V - \wt M$.
The simplex $S$ is extensible to a simplex $S'=(V_{S'},E_{S'})$
(with $S$ a sub-graph of $S'$) so that $V_{S'} = V_S \cup \{ a \}$ if
and only if there exists $\{ b,c \} \subseteq \wt M$ and two colors
$\wh A$ and $\wh B$ of $E - M$ so that $a b \in \wh A$ and $a c \in
\wh B$.

\end{corollary}

{\it Proof.} Use [C.3.6] and [T.3.9].

\section{Transitive orientations of an undirected graph}

\noindent In this section, using the results of the previous section
we prove the existence of a partition of maximal multiplices for the
set of edges of an undirected graph. Therefore, the transitive
orientations of a comparability graph turn up to the transitive
orientations of their multiplices. These orientations are independent
to each other. A theorem of decomposability for a non limit
undirected graph is proved.

\begin{lemma}\label{lem4.1}
Let $G=(V,E)$ be any graph. Let $X$ and $Y$ be two partitive sets of
$G$ so that $X \subseteq Y$. The following statements hold:

\begin{itemize}

\item[(i)] $X \in I_F (G) \Rightarrow X \in I_F (G (Y))\ ;$

\item[(ii)] $Y \in I_F (G) \Rightarrow (X \in I_F (G) \Leftrightarrow
X
\in I_F (G (Y)).$

\end{itemize}

\end{lemma}

{\it Proof.} Let $X$ be a ``strong'' partitive set of $G$ $(X \in I_F
(G))$. We have for every $Z \in I (G (Y)),\ Z \in I (G)$. Hence $X \in
I_F (G (Y))$. Let $X \in I_F (G (Y)),\ Y \in I_F (G)$ and $Z \in I
(G)$ so that $Z \cap X \not= \phi$ and $Z - X \not= \phi$. Then $Z
\cap X \not=
\phi \Rightarrow Z \cap Y \not= \phi$. But $Y \in I_F (G)$, thus
either
$Z \subseteq Y$ or $Y \subseteq Z$. But $(Y \subseteq Z \Rightarrow X
\subseteq Z)$ and $(X \subseteq Y \Rightarrow Z \in I (G (Y))
\Rightarrow X \subseteq Z)$.

\begin{lemma}\label{lem4.2}
Let $G=(V,E)$ be any graph. Let $F$ and $F'$ be two partitions of $G$
constituted of maximal ``strong'' partitive sets. We have then $F =
F'$.

\end{lemma}

{\it Proof.} Let us assume that $F \not= F'$. Since $F$ and $F'$ are
patitions of $V$, if $X \in F$ then there exists $X' \in F'$ so that
$X
\cap X' \not= \phi$. Thus either $X \subseteq X'$ or $X' \subseteq X$.
But
$X$ and $X'$ are maximal ``strong'' partitive sets, thus $X = X'$.
Hence $F = F'$.

\begin{proposition} \label{prop4.3}
Let $M (S)$ a multiplex of $G = (V,E)$ generated by the
simplex $S$ with rank$(M) \ge 2$. We have then for every $x \in \wt M$
there exists two colors $\{ \wh A,\wh B \} \subseteq \wh M$ with $\wt
A - \wt B \not= \phi$ and $\wt B - \wt A \not= \phi$ so that $x \in
\wt A
\cap
\wt B$.

\end{proposition}

{\it Proof.} Let $x \in M$, then there exists $\wh A \subset M$ and $y
\in \wt A$ so that $x y \in \wh A$. Moreover ($\wh A \subset M$ and
rank $(M) \ge 2$) $\Rightarrow S$ contains one tricolor triangle:
$(a, b, c\ ;\ \wh A, \wh B, \wh C)$ so that $b c \in \wh A,\ a c\in
\wh B$ and $a b \in \wh C$. Thus $a \not\in \wt A,\ b \not\in \wt B$
and $c
\not\in \wt C$~[Triangle lemma]. Hence $\wt A - \wt B \not= \phi$
and $\wt B - \wt A \not= \phi$. In the other hand:
$x y \in \wh A \Rightarrow ((a x \in \wh B\ \hbox{and}\ a y \in \wh
C)\
\hbox{or}\ (a x \in \wh C\ \hbox{and}\ a y \in \wh B))$~[Triangle
lemma]. If one suppose $a x \in \wh B$ then $x \in \wt A \cap \wt B$.

\begin{theorem}\label{th4.4}
Let $M (S)$ be a multiplex of $G = (V,E)$ generated by the simplex
$S=(V_S,E_S)$.
We have then the following statements:

\begin{itemize}

\item[(i)] $G (\wt M)$ has a partition of maximal ``strong''
partitive sets $F_M \not= \{ \wt M \}$~;

\item[(ii)] If rank$(M) = 1$, we have either $G (\wt M) / F_M$ is
isomorphic to $S$ or $G (\wt M) / F_M$ is indecomposable and
isomorphic to a sub-graph $G'=(V',E')$ of $G (\wt M)$ so that $E'
\subseteq M$~;

\item[(iii)] If rank $(M) \ge 2$, then $G (\wt M) / F_M$ is
isomorphic to $S$.

\end{itemize}

\end{theorem}

{\it Proof.} 1) If rank$(M) = 1$ then $M$ contains only one color $\wh
A = M$. Moreover if $G (\wt M)$ is indecomposable then $F_M = \{ \{ x
\}\ ;\ x \in \wt M \}$. Hence $G (\wt M) / F_M$ is isomorphic to $G
(\wt M)$ and thus $G (\wt M) / F_M$ is indecomposable.

Let us assume now that $G (\wt M)$ is decomposable. Thus $G
(\wt M)$ contains a non-trivial partitive set $X$.
After~[L.\rf{lem4.2}
$F_M$ exists. If $| F_M | = 2$ then $G (\wt M) / F_M$ is complete and
isomorphic to $S$.

Let us suppose that $| F_M | \ge 3$. If $G (\wt M) / F_M$ has a
non-trivial partitive set $X$, $G (\wt M) / F_M$ has a non-trivial
maximal ``strong'' partitive set $Y$~[L.\rf{lem4.2} and [C.3.5].

Hence $\cup Y \in I_F^{\star} (G (M))$[P.21] which is absurd because
$F_M$ is already maximal. Thus $G (\wt M) / F_M$ is indecomposable.
\vskip 0.8cm
2) Let us now assume that rank$(M) \ge 2$. Using the
proposition~[\rf{prop4.3} one gets that for every $x \in \wt M$ there
exists two colors $\{ \wh A, \wh B \} \subseteq \wh M$ with $\wt A -
\wt B \not= \phi$ and $\wt B - \wt A \not= \phi$ so that $x \in \wt A
\cap \wt B$.

Applying~[T.\rf{th3.3} we have that the intersection $\wt A \cap \wt B
\not= \wt M$ is a ``strong'' partitive set of $G$. Thus $\wt A \cap \wt
B$ is a ``strong'' partitive set of $G (\wt M)$ itself since $\wt M
\in I (G)$~[L.\rf{lem4.1}. Let $F_M$ be  the set of the intersections
two by two of colors of $M$.
Then for every $a \in V_S$ there exists $X \in F_M$ so that $a
\in X$. Let $\{ a, b \} \subset V_S$ so that $a \in X \in F_M$ and $b
\in Y \in F_M$. Let us assume that there exists a ``strong'' partitive
set $Z$ of $G (\wt M)$ so that $X \cup Y \subseteq Z$. $S$ will
contain a tricolor triangle $(a, b, c\ ;\ \wh A, \wh B, \wh C)$ with
$b c \in \wh A$, $a c \in \wh B$ and $a b \in \wh C$. This is absurd
since $\{ c \}$ and $Z$ are two ``strong'' partitive sets of $G (\wt
M)$ and can be related at most by only one color[T.3.4]. Hence $F_M$
is a maximal partition of "strong" partitive sets and
separates the summits of $S$.
 Finally we get that $G
(\wt M) / F_M$ is isomorphic to $S$.

\begin{corollary}\label{cor4.5}
The only multiplices $M (S)$ which might be  not transitively
orientable are those of
rank $= 1$ and so that $G (\wt M) / F_M$ is non isomorphic to $S$.

\end{corollary}

{\it Proof.} Because complete graphs are orientable.

\begin{corollary}\label{cor4.6}
An undirected graph $G = (V, E)$ can have at most one multiplex which
spanned all its summits.

\end{corollary}

{\it Proof.} Let $M (S)$ and $M (S')$ be two multiplices of $G$ so
that $\wt M = V$ and $\wt M' = V$, we have then $G (\wt M) = G (\wt
M') = G$. Let $F$ and $F'$ be two partitions of maximal ``strong''
partitive sets of $G$ related to $M$ and $M'$. Thus $G / F$ is
isomorphic to $S$ and $G / F'$ is isomorphic to $S'$. Moreover the two
isomorphismes conserve the colors~[C.\rf{cor3.5}. But
after~[L.\rf{lem4.2} we have $F = F'$. Hence $S = S'$ and $M = M'$.

\begin{corollary}\label{cor4.7}
Let $G = (V,E)$ be an undirected graph and $M (S),\ M' (S')$ two
maximal multiplices of $G$. We have then $M \cap M' \not= \phi
\Rightarrow M = M'$.

\end{corollary}

{\it Proof.} Since $M$ and $M'$ are ``strong'' partitive sets we
have \\
$M \cap M' \not= \phi \Rightarrow \wt M\cap \wt {M'}\neq \phi
\Rightarrow ( \wt M \subseteq \wt M'\ \hbox{or}\ \wt M' \subseteq \wt
M )$.\newline Let us assume that $\wt M \subseteq \wt M'$ and let $F,
F'$ be two partitions of maximal ``strong'' partitive sets related
respectively to
$G (\wt M)$ and $G (\wt M')$. If there exists $X \in F'$ so that $M
\cap E (X) \not= \phi$ then $M \subseteq E (X)$~[T.\rf{th3.7} and $M
\cap M' = \phi$ ( since $G (M') / F'$ is isomorphic to $S'$ and $S$
would be a sub-graph of $G (X)$. Hence $E_S\cap E_{S'}=
\phi$).\newline Finaly we have that for every $X \in F',\ E (X) \cap M
= \phi$. Thus
$S$ is isomorphic to a sub-graph of $S'$. But $S$ is maximal. Hence $S
= S'$ and $M = M'$.

\begin{corollary}\label{cor4.8}
Let $G = (V,E)$ be an undirected graph with $E \not= \phi$. $E$ has
then a partition of maximal multiplices.

\end{corollary}

The theorem 4.4 tell us that a multiplex has the same number of
transitive orientations as the simplex which generated
this multiplex. The simplex itself has a number of transitive
orientations equal to the number of the possible permutations of
its summits. Moreover [C.\rf{cor4.5} asserts that the only multiplices
$M(S)$ which might  be not transitively orientable are those with rank
1 so that $G(\wt M)/F_M$ is not isomorphic to $S$. Thus using
[C.\rf{cor4.8} the problem of transitive orientation for a
comparability graph come down to the transitive orientation of its
multiplices. But the following problem is rised: if we orientate in
any way and at certain step a given multiplex, will this orientation
influence or not the orientations of the other multiplices at the
following steps ? The response is not and we will prove this statement
using a theorem [T.\rf{th4.11} which is known and for which we propose
a new proof outcoming {}from the forcing theory.\newline
Before announcing [T.\rf{th4.11}, we will announce a theorem
[T.\rf{th4.9} which, in fact, is a mathematical algorithm permitting
to find all the transitive orientations for a comparability graph
which has only non limit sub-graphs, {\em e.g.}, case of finite
graphs.

\begin{lemma}\label{lem4.x}
Let $G=(V,E)$ be  a connected undirected graph. Then $G$ can not
contain a multiplex $M$ so that both $\wt M\neq V$ and $\wt M$ is
maximal for the inclusion among the $\wt N$, where $N$ is any
multiplex of $G$

\end{lemma}

{\it Proof.} Let us assume that such a multiplex $M$ exists and
show that it is absurd. Since $\wt M$ is maximal for the inclusion it
implies that $M$ is maximal. After[T.\rf{th3.9}, $\wt M$ is a "strong"
partitive set. Since $\wt M\neq V$  and $G$ is connected, we have
:\newline there exists $x\in V-\wt M$ and $y\in \wt M$ so that $xy\in
E$. Let $\wh A$ the color containing $xy$. Then
$\wh A\not\subset M$ and $x\in (\wt A \cap \wt M )\neq\phi$.
But $\wt M$ is a "strong" partitive set of $G$ and $\wt A$ is a
partitive set of $G$, thus $\wt M\subset \wt A$, which is absurd
because $\wh A$ is a multiplex of rank 1.

\begin{theorem}\label{th4.9}
Let $G = (V,E)$ be an undirected graph having a partition of maximal
``strong'' partitive sets $F_G \not= \{ V \}$. We have then that $G /
F_G$ satisfies one of the following exclusive assertions:

\begin{itemize}

\item[(i)] $G / F_G$ is empty.

\item[(ii)] $G / F_G$ is indecomposable and there exists a maximal
multiplex $M_G$ of $G$ with  rank$(M_G) = 1,\ \wt M_G = V$ and $G /
F_G$ is isomorphic to a sub-graph $G'=(V',E')$ of $G$
so that $E' \subseteq
M_G$.

\item[(iii)] $G / F_G$ is complete and isomorphic to a maximal simplex
$S$ generating a maximal multiplex $M_G$ so that $\wt M_G = V$.

\end{itemize}

\end{theorem}

{\it Proof.} If $G$ is non connected, $F_G$ is the class of the
connected components and $G / F_G$ is empty. In the other hand, it is
obvious that if $G / F_G $ is empty then $G$ is non connected. Hence
$G / F_G$ is empty if and only if $G$ is non connected.

Let us assume now that $G$ is connected. Since $G$ is non limit,
it implies that $G$ has maximal multiplex $M$ so that $\wt M$
is maximal for the inclusion. Thus using the [L.\rf{lem4.x} we have
$\wt M =V$. Hence $G=G(\wt M)$ and using [T.\rf{th4.4} we get the
result.

\begin{corollary}\label{cor4.10}
Let $G = (V,E)$ be a connected and undirected graph. $G$ has then a
partition of maximal ``strong'' partitive sets $F_G \not= \{ V \}$ if
and only if $G$ has a multiplex $M_G$ so that $\wt M_G = V$.

\end{corollary}

{\it Proof.} Applying~[T.\rf{th4.9} we have: $F_G$ exists $\Rightarrow
M_G$ exists. In the other hand applying~[T.\rf{th4.4}, we have the
other implication: $M_G$ exists $\Rightarrow F_G$ exists.

\begin{theorem}\label{th4.11}
Let $O = (V,E)$ be a partial order and $G = (V, E)$ be its
comparability graph. $O$ and $G$ have then the same ``strong''
partitive sets $(I_F (G) = I_F (O))$.

\end{theorem}

Before giving the proof of this theorem, we present some preleminary
results which will be used for the proof.

\begin{lemma}\label{lem4.12}
Let $O = (V,E')$ be a partial order and $G = (V, E)$ its
comparability graph. Then every partitive set of $O$ is a partitive
set of $G$\\
$ (I (O) \subseteq I (G))$.

\end{lemma}

{\it Proof.} Let $Y \in I (O)$. If $Y$ is a trivial partitive set, we
have $Y \in I (G)$. Let us suppose that $Y$ is non trivial, then if $a
\in V - Y$ we have one of the following statements:

\begin{itemize}

\item for every $y \in Y,\ ((a,y) \in E')\Rightarrow a y \in E$.

\item or for every $y \in Y,\ ((y,a) \in E')\Rightarrow a y \in E$.

\item or for every $y \in Y,\ (\{ (a,y), (y,a) \} \cap E' = \phi)
\Rightarrow a y \not\in E$.

\end{itemize}

Therefore $Y \in I (G)$.

\begin{proposition}\label{prop4.13}
Let $O = (V,E')$ be a partial order and $G = (V, E)$ its
comparability graph. Then every ``strong'' partitive set of $G$ is a
``strong'' partitive set of $O$  $ (I_F (G) \subseteq I_F (O))$.

\end{proposition}

{\it Proof.} Let $X$ be a ``strong'' partitive set, {\it i.e.,} $X \in I_F
(G)$. If $X$ is trivial then $X \in I_F (O)$. Let us
suppose that $X$ is non trivial. We have then one of the following
statements:

\begin{itemize}

\item for every $x \in X,\ x a \not\in E$, then for every $\{ x, y \}
\subseteq X,\ (x, a) \cong (y,a)$ (cf.~\S2.1).

\item or there exists a color $\wh A$ of $G$ so that for every $x \in
X,\ x a \in \wh A$, then for every $\{ x, y \} \subseteq X,\ (x, a)
\cong (y, a)$ (since $G$ is a comparability graph [T.\rf{th2.2}).
\\ Therefore $X \in I (O)$. Let $Y \in I (O)$ so that $Y \cap X \not=
\phi$ and $Y - X \not= \phi$. Then using the previous
lemma~[\rf{lem4.12} we have $Y \in I (G)$ and then $X \subset Y$. Thus
$X \in I_F (O)$.

\end{itemize}

\begin{proposition}\label{prop4.14}
If $X$ is a maximal ``strong'' partitive set of $G$ then $X$ is a
maximal ``strong'' partitive set of $O$.

\end{proposition}

{\it Proof.} If $X \in I_F (G)$ is maximal, it implies that $G$ is non
limit. Thus $G$ has a partition ${ F}$ of maximal ``strong''
partitive sets and $X \in { F} \subseteq I_F (O)$.
After~[T.\rf{th4.9} we have either $G / { F}$ is empty and then $X$
is a connecting class. Thus $X$ is maximal in $I_F (O)$ ; or $G / {
F}$ is complete and then $O / {  F}$ is a chain, therefore $G /
{  F}$ do not has a non trivial ``strong'' partitive sets, thus $X$ is
maximal in $I_F (O)$ ; or $G / {  F}$ is indecomposable and the
partitive sets of $G / {  F}$ are trivial, therefore,
after~[L.\rf{lem4.12} the partitive sets of $O / {  F}$ are trivial,
thus $X$ is maximal in $I_F (O)$.

\begin{proposition}\label{prop4.15}
Let $X$ be a non trivial ``strong'' partitive set of $O$. There exists
then a non trivial ``strong'' partitive set $Y$ of $G$ so that
$X \subseteq Y$.

\end{proposition}

{\it Proof.} If $G$ is non limit then it has a partition ${  F}$ of
maximal ``strong'' partitive sets which is also a partition of maximal
``strong'' partitive sets of $O$[P.4.15]. Thus there exists $Y \in {
F}$ so that $X \subseteq Y$. If $G$ is  limit then since $X \in I (O)
\Rightarrow X \in I (G)$~[L.\rf{lem4.12}, therefore $Y$ exists.

\vskip 0.8cm
\noindent
In what follows, we present the proof of the theorem~[\rf{th4.11}.

{\it Proof.} If $X \in I_F (G)$ then after~[P.\rf{prop4.14} $X \in I_F
(O)$.
Let $X \in I_F^{\star} (O)$. Then after~[P.\rf{prop4.15} there exists
$Y
\in I_F^{\star} (G)$ so that $X \subseteq Y$. Let us suppose that $X
\not\in I_F (G)$. Thus $X \subset Y$.\newline
$Y \in I_F^{\star} (G) \Rightarrow Y \in I_F^{\star} (O)$. Therefore
after~[L.\rf{lem4.1} $X \in I_F^{\star} (O (Y))$. Thus there exists
$Y_1
\in I_F^{\star} (G (Y))$ so that $X \subset Y_1$. So we have
constructed a strictly decreasing suite of elements of $I_F^{\star}(G)
\subseteq I_F^{\star} (O)$ (after~[P.\rf{prop4.13}).\newline
The intersection $\cap Y_i$ of this suite is a ``strong'' partitive
set of $O$ and it is  the smallest element of $I_F^{\star} (O)$ which
contains
$X$. Therefore $X = \cap Y_i$. But $\cap Y_i$ is also a ``strong''
partitive set of $G$. Thus $X \in I_F^{\star} (G)$ and then $I_F (O) =
I_F (G)$.

\section*{Acknowledgments}
One of the authors, M.~H., would like
to thank Prof. P. ILLE for helpful and stimulating discussions
and also Profs. C. Rauzy and G. Fardoux for their encouragements.

\end{document}